\documentclass[a4paper,11pt,UKenglish]{article}

\usepackage{amsmath,amssymb, xcolor, listings, amsthm}
\usepackage{mathtools}
\usepackage{enumitem}
\usepackage{url}
\usepackage{hyperref}
\usepackage{fullpage}
\usepackage{diagbox}
\usepackage{tikz}

\usepackage{graphicx}
\usepackage{caption}
\usepackage{indentfirst}
\usepackage{mathrsfs}
\usepackage{manfnt}
\usepackage[capitalize]{cleveref}

\usepackage{todonotes}

\usepackage[normalem]{ulem}

\hypersetup{colorlinks=true, linkcolor=black, citecolor=blue, urlcolor=blue}
\theoremstyle{definition}
\newtheorem{definition}{Definition}[section]
\newtheorem{oss}[definition]{Remark}
\newtheorem{example}[definition]{Example}

\theoremstyle{plain}
\newtheorem{theorem}[definition]{Theorem}
\newtheorem{prop}[definition]{Proposition}

\newtheorem{lemma}[definition]{Lemma}
\newtheorem{cor}[definition]{Corollary}
\newtheorem{remark}[definition]{Remark}

\DeclareMathOperator{\bub}{\mathbf{B}}

\title{Preimages under the bubblesort operator}
\author{Mathilde Bouvel, Lapo Cioni, Luca Ferrari}
\date{}

\begin{document}

\maketitle

\begin{abstract}
We study preimages of permutations under the bubblesort operator $\mathbf{B}$. We achieve a description of these preimages much more complete than what is known for the more complicated sorting operators $\mathbf{S}$ (stacksort) and $\mathbf{Q}$ (queuesort). 
We describe explicitly the set of preimages under $\mathbf{B}$ of any permutation $\pi$ from the left-to-right maxima of $\pi$, showing that there are $2^{k-1}$ such preimages if $k$ is the number of these left-to-right maxima. 
We further consider, for each $n$, the tree $T_n$ recording all permutations of size $n$ in its nodes, in which an edge from child to parent corresponds to an application of $\mathbf{B}$ (the root being the identity permutation), and we present several properties of these trees. In particular, for each permutation $\pi$, we show how the subtree of $T_n$ rooted at $\pi$ is determined by the number of left-to-right maxima of $\pi$ and the length of the longest suffix of left-to-right maxima of $\pi$. 
Building on this result, we determine the number of nodes and leaves at every height in such trees, and we recover (resp. obtain) the average height of nodes (resp. leaves) in $T_n$. 
\end{abstract}

\section{Introduction}

\subsection{Motivation}

The foundational work of Knuth~\cite[Section 2.2.1]{knuth1} defines the stacksorting procedure and relates it to pattern-avoidance in permutations. 
Since then, many similar sorting procedures have been defined, where the sorting device is not necessarily a stack (\emph{e.g.}, a queue, a deque, \dots) and where such devices can be combined. 
These sorting procedures have been studied from various points of view, and we refer the reader to the surveys~\cite{BonaSurvey,LucaTalk} and to the introduction of~\cite{VatterSurvey} for an overview of the research on this topic and bibliographic pointers. 

Our focus in the present work is on the bubblesort operator, $\bub$, which corresponds to one pass of the bubblesort algorithm. 
Although not strictly speaking a sorting procedure associated with a particular device, $\bub$ shares some features with such procedures, as demonstrated in~\cite{bubbleM}.  
Our point of focus is the study of preimages of permutations under $\bub$, 
a topic which has proved very rich for other sorting operators. We start by reviewing the literature related to the study of preimages under a sorting operator, focusing on different aspects of this study. 

\paragraph{Existence and number of preimages}
To our knowledge, Bousquet-Mélou is the first to investigate a sorting operator (here, the stacksorting operator $\mathbf{S}$) through the lens of preimages. In her article~\cite{Mireille}, she provides an algorithm to decide if a permutation has at least one preimage under $\mathbf{S}$, and in this case to compute one such preimage called \emph{canonical}, from which all other preimages can then be computed. She also derives an equation for the generating functions of permutations having at least one preimage under $\mathbf{S}$.

Recent improvements on the study of permutations having preimages under $\mathbf{S}$ have been obtained by Defant and co-authors. 
Specifically, \cite{defant2} counts permutations having exactly one preimage under $\mathbf{S}$; \cite{defant3} studies the number of preimages that a permutation can have; and~\cite{defant4} uses a method involving the computation of preimages to give bounds on the number of $t$-stack-sortable permutations of size $n$. 

In~\cite{LucaLapo1,LucaLapo2}, Cioni and Ferrari consider another sorting operator: the queuesort operator $\mathbf{Q}$. They provide a recursive description of the preimages under $\mathbf{Q}$, study the possible numbers of preimages that a permutation can have, and compute the number of permutations having 0, 1 or 2 preimages. 

\paragraph{Sorting trees}
In~\cite{Mireille}, Bousquet-Mélou also defined the \emph{sorting trees} associated with $\mathbf{S}$. For each integer $n$, the sorting tree (see~\cite[fig. 2]{Mireille}) for size $n$ is the tree whose nodes are the permutations of size $n$, with root the identity permutation $1 2 \dots n$, and such that the children $\tau$ of any permutation $\sigma$ are those such that $\mathbf{S}(\tau) = \sigma$. The article~\cite{Mireille} already describes some properties of the sorting trees (although not phrased as such, the focus in~\cite{Mireille} being different). 

Building on these, Defant~\cite{defant7} proved in disguise further properties of the sorting trees. Namely, assigning to any permutation $\sigma$ a ``label'' which consists of the skeleton of the decreasing binary tree whose in-order reading is $\sigma$, he shows~\cite[Theorem 11]{defant7} that the label of $\sigma$ determines the number and the labels of the children of $\sigma$ in the sorting tree, and hence recursively the whole skeleton of the subtree of the sorting tree rooted at $\sigma$. 

We point out an easy consequence of this fact: for every permutation $\sigma$, there exists a permutation $\tau$ which is at distance one from the root in the sorting tree, and such that the subtrees of the sorting tree rooted at $\sigma$ and at $\tau$ are isomorphic. This indeed follows from the two facts that permutations at distance one from the root in the sorting tree are the stack-sortable permutations (\emph{i.e.}, those avoiding $231$) and that every binary tree has a decreasing labeling which avoids $231$. 

About queuesort,~\cite{LucaLapo1} mentions that it would be interesting to study the properties of the sorting trees associated with $\mathbf{Q}$, but we are not aware of any such results at the moment. 

\paragraph{Preimages of permutation classes} 
The permutation class defined by a set of excluded patterns $B$ is the set of all permutations avoiding every pattern in $B$, and is denoted $\textnormal{Av}(B)$. A consequence of the characterization of stack-sortable permutations as those avoiding $231$ is that the permutations sortable by two applications of $\mathbf{S}$ (called West-two-stack-sortable) are the preimages of permutations in the class $\textnormal{Av}(231)$.
Describing the West-two-stack-sortable permutations is therefore an instance of the more general question asking for a description of the preimages (for $\mathbf{S}$) of a class of pattern-avoiding permutations. 

This question has been studied algorithmically in~\cite{ClaessonUlfarsson}, where Claesson and Ulfarsson provide an algorithm to describe preimages of principal pattern classes, in terms of decorated patterns. 
Their approach has been extended in~\cite{Magnusson} to other sorting operators, including $\mathbf{Q}$. 

A more enumerative perspective on preimages of permutation classes for $\mathbf{S}$ has later been provided in~\cite{defant5,defant6}. 

Finally, the problem of describing preimages of permutation classes for $\bub$ has been completely solved in~\cite{bubbleM} for principal permutation classes. 

\paragraph{Complexity of sorting procedures}
As sorting algorithms, the procedures considered above are clearly inefficient. Nevertheless, it makes sense to ask how many applications of a sorting operator (like $\bub$ or $\mathbf{S}$ or others) are needed to fully sort a permutation of size $n$, either in the worst case or on average. Some bounds for the average case are provided in~\cite{Tao} for various sorting operators, and in~\cite{defant1} for 
$\mathbf{S}$. 

\subsection{Our results}

We focus on the bubblesort operator $\bub$, and study preimages of permutations under $\bub$. As previously indicated, a study of preimages of permutation classes has already been done in~\cite{bubbleM}, and we leave these aspects aside of our study. 

About the first point above (characterizing the existence and number of preimages), we can be very precise when confronting these questions with bubblesort. Indeed, we can fully describe the set $\bub^{-1}(\sigma)$ of preimages of any permutation $\sigma$. This is presented in \cref{sec:CompPreimages}. 

The description of $\bub^{-1}(\sigma)$ is actually rather simple, and involves essentially the left-to-right maxima of $\sigma$. 
From this description, we deduce that the number of preimages of a permutation $\sigma$ is $2^{k-1}$ for $k$ the number of left-to-right maxima of $\sigma$ if $\sigma$ ends with its maximum (and $0$ otherwise).  
The results are therefore much more precise than what can be achieved with $\mathbf{S}$ for instance. An informal explanation which we can offer to explain this fact is the following. The operator $\bub^{-1}$ acts only on the left-to-right maxima of a permutation, and does so in a very controlled way, allowing to describe the left-to-right maxima after the application of $\bub^{-1}$. This is very useful in particular for describing iterated preimages for $\bub$. 

In \cref{sec:Tree}, we turn to the study of these iterated preimages. Specifically, we define for each permutation $\pi$ the tree $T(\pi)$ whose root is $\pi$ such that the children of any permutation are the preimages of this permutation. These trees are the analogue in the case of $\bub$ of the \emph{sorting trees} defined by Bousquet-Mélou in~\cite{Mireille} for $\mathbf{S}$.

As expected from the informal discussion above, we show that these trees are completely determined by the left-to-right maxima of $\pi$. More precisely, they are determined by what we call the \emph{label} of $\pi$, which plays the same role as the ``label'' we presented above in the case of $\mathbf{S}$ but is much simpler than it (recall that this was the skeleton of the decreasing binary tree whose in-order reading is $\pi$). In the case of $\bub$, the \emph{label} of $\pi$ is the pair consisting of the number of left-to-right maxima of $\pi$ and the length of the longest suffix of left-to-right maxima of $\pi$. 
We also prove, similarly to the case of $\mathbf{S}$, that every tree $T(\pi)$ is isomorphic to a tree $T(\tau)$ for $\tau$ a permutation such that $\bub(\tau)$ is the identity.

Finally, in \cref{sec:height} we study heights of nodes and leaves in $T(\pi)$. For $\pi$ the identity permutation of size $n$, the average height of a node in $T(\pi)$ is the average number of passes of $\bub$ necessary to sort a permutation of size $n$, which is known from~\cite[Theorem 7.14]{FS} for instance. We modify this analysis to also compute the average height of a leaf in $T(\pi)$ (which corresponds to the number of passes of $\bub$ necessary to sort a permutation of size $n$ which does not belong to the image of $\bub$). 
In addition, for any permutation $\pi$, we provide closed formulas for the number of nodes and leaves at any possible height $j$ in $T(\pi)$; these formulas depend only on $j$ and on the label of $\pi$. 
Nevertheless, we could not deduce the average height of a leaf or node in $T(\pi)$ from these formulas. 

Before moving to definitions and basic properties of $\bub$, we note that we see the questions studied here as bubblesort analogues of similar questions previously studied on the stacksort and queuesort operators. The operator $\bub$ being simpler, the answers obtained are much more precise than in the case of $\mathbf{S}$ or $\mathbf{Q}$. 
A possible direction for future research is to confront these questions with other sorting operators as well. Those listed in~\cite{Magnusson} can be a very good source of inspiration.

\subsection{The operator $\bub$: definition and some basic properties}

For any integer $n$, a permutation $\pi$ of size $n$ is a sequence $\pi_1 \pi_2 \dots \pi_n$ containing exactly once each symbol from $\{1,2, \dots, n\}$. 
An element $\pi_j$ of a permutation $\pi$ is a left-to-right maximum if it is larger than all elements to its left, that is to say $\pi_j$ is such that $\pi_j > \pi_i$ for all $i<j$. 

The bubblesort operator, denoted $\bub$, corresponds to applying one pass of the classical bubblesort algorithm to a permutation. Specifically, $\bub(\pi)$ is obtained from $\pi$ scanning its elements from left to right, each time exchanging an element with the one sitting to its right whenever the latter is smaller. Thereby (see also \cref{lem:image} below), the left-to-right maxima of the permutation ``bubble up'' to the right, until they are blocked by the next left-to-right maximum. 

For example, for $\pi = {\bf 4} 2 1 {\bf 6} 3 {\bf 7} {\bf 8} 5$, the left-to-right maxima are $4,6,7$ and $8$ (shown in bold) and $\bub(\pi) = 2 1 4 3 6 7 5 8$.

\begin{oss}
\label{rk:not_perm}
From the above definition, it is clear that $\bub$ can be applied \emph{verbatim} to sequences of distinct integers which are not necessarily permutations. Consequently, every statement about $\bub$ on permutations also applies to sequences of distinct integers up to relabeling the values with the order-isomorphic permutation. 

For instance, $\bub({\bf 4} 2 1 {\bf 6} 3) = 2 1 4 3 6$, and the preimages of $2 1 4 3 6$ for $\bub$ are obtained from the preimages of the permutation $21435$ replacing $5$ with $6$. 
\end{oss}

The bubblesort operator may be described in several other ways, and we give two below. The reader needing an explanation of the equivalence with the above definition can find it in~\cite[Lemma 1 and the equation displayed just above it]{bubbleM}. 
The first parallels the recursive definition of the stacksorting operator $\mathbf{S}$: decomposing $\pi$ into $\pi = \pi_L n \pi_R$ with $n$ the maximal value occurring in $\pi$, we have $\bub(\pi) = \bub(\pi_L) \pi_R n$ (and by comparison $\mathbf{S}(\pi) = \mathbf{S}(\pi_L)\mathbf{S}(\pi_R) n$).
The second focuses on the left-to-right maxima, and we record it in a lemma for future reference. 

\begin{lemma} \label{lem:image}
Let $\pi$ be a permutation, and write $\pi = \mu_1 A_1 \mu_2 A_2 \dots \mu_k A_k$, 
where the $\mu_i$'s are all the left-to-right maxima of $\pi$ (and the $A_i$ are possibly empty sequences of integers). 
Then $\bub(\pi) = A_1 \mu_1 A_2 \mu_2 \dots A_k \mu_k$.
\end{lemma}

In particular, a permutation is in the image of $\bub$ if and only if it ends with its maximum. 
We record below another observation which follows immediately from~\cref{lem:image}

\begin{cor}\label{cor:subsetOfLRMax}
For any $\pi$, 
the set of left-to-right maxima of $\pi$ is included in the set of left-to-right maxima of $\bub(\pi)$. 
\end{cor}

\section{Computing the preimages} 
\label{sec:CompPreimages}

In this section we present a procedure to compute the set $\bub^{-1}(\sigma)$ of all preimages of any given permutation $\sigma$. 
First, as noted just after \cref{lem:image}, 
for a permutation $\sigma$ which does not end with its maximum, $\bub^{-1}(\sigma)$ is empty. This trivial case being solved, we now focus on the interesting case where $\sigma$ does end with its maximum. 

Let $\sigma = \sigma_1 \sigma_2 \dots \sigma_n$ be a permutation of size $n$ which ends with its maximum. 
Define $P$ as a set which contains only $\sigma$. For each $i$ from $n$ down to $2$, do the following: for each $\pi \in P$, 
\begin{itemize}
 \item if $\pi_{i-1}$ is not a left-to-right maximum of $\pi$, then replace $\pi$ in the set $P$ by the permutation $\pi_1 \dots \pi_i \pi_{i-1} \dots \pi_n$ (that is to say, we swap $\pi_i$ and $\pi_{i-1}$); 
 \item if $\pi_{i-1}$ is a left-to-right maximum of $\pi$, then $\pi$ stays in the set $P$, and in addition we add in $P$ the permutation $\pi_1 \dots \pi_i \pi_{i-1} \dots \pi_n$ (where $\pi_i$ and $\pi_{i-1}$ are swapped).
\end{itemize}

\begin{example}
The table below shows the evolution of the set $P$ of the above procedure, for $\sigma = 325146$. 
\begin{center}
\begin{tabular}{|c|c|c|c|c|c|c|}
\hline
 $i= \dots$ & initialization & $6$ & $5$ & $4$ & $3$ & $2$ \\
\hline 
$P$ contains & $325146$ & $325164$ & $325614$ & $325614$ & $352614$ & $352614, 532614$\\ 
& & & & $326514$ & $362514$ & $362514, 632514$\\ 
\hline
\end{tabular}
\end{center}
\end{example}

\medskip

We may note that, when starting any step $i$, for any $\pi \in P$, 
$\pi_i$ is always a left-to-right maximum of $\pi$, and
$\pi_{i-1}$ is a left-to-right maximum of $\pi$ if and only if $\sigma_{i-1}$ is a left-to-right maximum of $\sigma$. 
Indeed, all steps until step $i$ (excluded) of the above procedure leave the prefixes of length $i-1$ unchanged. 

It is useful to have a different (although equivalent) 
presentation of this procedure, which we now give. 
Starting from $\sigma$, where we see the rightmost element $\sigma_n$ as distinguished, 
we move the distinguished element to the left until it becomes the leftmost, 
according to the following rules.
\begin{itemize}
 \item If the element immediately to the left of the distinguished one is not a left-to-right maximum of the current sequence, 
 then the distinguished element is forced to move to the left (\emph{i.e.} is swapped with its left neighbor). The distinguished element remains the same. 
 \item If the element immediately to the left of the distinguished one is a left-to-right maximum of the current sequence, 
 then the distinguished element may either move to the left or stay in place. 
 In the first case, the distinguished element remains the same. 
 In the second case, the distinguished element becomes the left neighbor of the previously distinguished element. 
\end{itemize}
It is easy to see that the set $P$ computed by the original procedure consists of 
all possible results of applying this alternative procedure. 

Some remarks (all easily observed) about this alternative procedure are useful. 
First, the index $i$ (between $n$ and $2$) of a given step of the original procedure
always corresponds to the position of the distinguished element in the evolving sequence. 
Second, the distinguished element is always a left-to-right maximum of $\sigma$, and also of the evolving sequence. 
Third, for any sequence $\pi$ produced, the elements which were at some point distinguished in the sequence are exactly the left-to-right maxima of $\pi$.

\begin{theorem}\label{thm:SetOfPreimages}
Let $\sigma$ be any permutation ending with its maximum, and $P$ be the set produced by the procedure above. 
Then $P$ is the set of preimages of $\sigma$ under $\bub$, that is to say, $P=\bub^{-1}(\sigma)$.
\end{theorem}

\begin{proof}
Assume first that a sequence $\pi$ has been produced by the above procedure. 
It means that $\pi$ has been produced from $\sigma$ by considering some left-to-right maxima of $\sigma$, from the right to the left, and moving these left-to-right maxima to the left, until they reach the position of the next left-to-right maximum which moves to the left. 
This is exactly undoing the action of $\bub$. 
More precisely, assume that $\pi=\mu_1 A_1\cdots \mu_k A_k$ with the $\mu_i$ the left-to-right maxima of $\pi$. 
By construction of $\pi$, $A_i$ is not empty if and only if $\mu_i$ has been moved by our procedure. 
Applying $\bub$ to $\pi$ yields $ A_1 \mu_1 \cdots A_k \mu_k$, 
thus exactly the elements that were moved by our procedure will be moved to the right by bubblesort. 
Moreover, $\bub$ moves $\mu_i$ to the right until it reaches the position immediately before $\mu_{i+1}$, 
and we claim that this is $\mu_i$'s original position in $\sigma$: 
indeed, our procedure only moves the $\mu_j$'s, so it must have started moving $\mu_i$ to the left immediately after considering $\mu_{i+1}$.
This proves that $\bub(\pi) = \sigma$, 
and therefore $P \subseteq \bub^{-1}(\sigma)$. 

For the converse inclusion, we proceed by induction on the size of $\sigma$. 
The statement is obvious for size $1$. 
So, let us consider $\pi \in \bub^{-1}(\sigma)$, for $\sigma$ of size greater than $1$.
We decompose $\pi$ around its maximal element as $\pi = L n R$. 
Then $\bub(\pi) = \bub(L) R n$, so that $\sigma = \bub(L) R n$. 
Starting from $\sigma = \bub(L) R n$, the above procedure is always allowed to move $n$ towards the left, 
and may decide when reaching $\bub(L) n R$ to distinguish the last element of $\bub(L)$ instead of $n$. 
Indeed, $\bub(L)$ necessarily ends with its maximum, 
so that the last element of $\bub(L)$ is a left-to-right maximum. 
Since $\bub(L)$ is a sequence shorter than $\sigma$ ending with its largest value, 
we may apply the induction hypothesis to it, and deduce that $L$ has been produced by the above procedure applied to $\bub(L)$.  
Combining these two facts, it follows that 
our procedure applied to $\sigma = \bub(L) R n$ can produce $L n R = \pi$,  
therefore showing that $\bub^{-1}(\sigma) \subseteq P$.
\end{proof}

\cref{thm:SetOfPreimages} has several consequences. 
First, we can refine \cref{cor:subsetOfLRMax} 
and describe $\bub^{-1}(\sigma)$ exactly from the left-to-right maxima of $\sigma$. 

\begin{cor}\label{cor:numberOfPreimages}
Let $\sigma$ be a permutation of size $n$ ending with its maximum (\emph{i.e.}, $\sigma_n=n$). 
Let $k$ be the number of left-to-right maxima of $\sigma$ (including $n$). 

There is a bijective correspondence between 
the preimages of $\sigma$ under $\bub$ and
the subsets of the $k-1$ left-to-right maxima of $\sigma$ different from $n$. 

More precisely, this correspondence works as follows. 
For any set $S = \{s_1 < \dots < s_j\}$ of $j \leq k-1$ left-to-right maxima of $\sigma$ different from $n$, 
writing $\sigma = B_0 s_1 B_1 s_2 B_2\dots s_{j-1} B_{j-1} s_j B_j n$ (for the $B_i$ possibly empty sequences of integers, 
which contain the $k-j$ left-to-right maxima not in $S$
and the elements of $\sigma$ which are not left-to-right maxima), 
the corresponding preimage of $\sigma$ is 
$ s_1 B_0 s_2 B_1 \dots s_j B_{j-1} n B_j$. 
\end{cor}

\begin{proof}
From the alternative description of the procedure computing $\bub^{-1}(\sigma)$, 
we have seen that 
the elements which are distinguished at some point are exactly the left-to-right maxima of the preimage produced. 
In addition, by definition of this procedure, 
the distinguished elements form a subset containing $n$ of the set of left-to-right maxima of $\sigma$. 
This proves the claimed bijective correspondence.

To describe precisely the preimage corresponding to a subset $S$, 
it is enough to note that every distinguished element moves to the left 
until a new distinguished element is chosen, leaving all other elements unchanged. 
\end{proof}

This allows to count the preimages of any given permutation, 
in total or by the number of their left-to-right maxima. 

\begin{cor}\label{cor:PreimagesFixedNumberLTR}
Let $\sigma $ be a permutation of size $n$ ending with its maximum, and with $k$ left-to-right maxima. 

The cardinality of $\bub^{-1}(\sigma)$ is $2^{k-1}$, 
and for any $1 \leq j \leq k$, 
the number of preimages of $\sigma$ with $j$ left-to right maxima is $\binom{k-1}{j-1}$.
\end{cor}
\begin{proof}
The cardinality of $\bub^{-1}(\sigma)$ follows immediately from \cref{cor:numberOfPreimages}. 
By \cref{cor:numberOfPreimages}, 
a preimage of $\sigma$ with $j$ left-to-right maxima corresponds bijectively to 
a subset containing $j-1$ elements of the set of left-to-right maxima of $\sigma$ different from $n$. We have $\binom{k-1}{j-1}$ different ways to select these subsets, thus proving the lemma.
\end{proof}

Second, we can characterize the permutations having a given number of preimages. 

\begin{cor}
For any $k \geq 1$, the permutations having exactly $2^{k-1}$ preimages under $\bub$ are those ending with their maximum and having $k$ left-to-right maxima in total. 

In particular, there are $\left[ {n-1} \atop {k-1} \right]$ permutations of size $n$ having $2^{k-1}$ preimages under $\bub$, where $\left[ {n} \atop {k} \right]$ are the (unsigned) Stirling numbers of the first kind. 
\end{cor}

\begin{proof}
The first statement follows immediately from \cref{cor:numberOfPreimages}. 
The second follows from the well-known fact that 
Stirling numbers of the first kind enumerate permutations according to their size and number of cycles, 
using the classical Foata bijection which maps 
permutations of size $n$ with $k$ cycles to permutations of size $n$ with $k$ left-to-right maxima. 
\end{proof}

\section{The trees of iterated preimages} 
\label{sec:Tree}

For any $n$, we denote by $S_n$ the set of permutations of size $n$ and by $id_n = 1 2 \dots n$ the identity permutation of size $n$. 
We start by defining $T(\pi)$ for any permutation $\pi$, and $T_n = T(id_n)$.

\begin{definition}
\label{def:tree}
Let $T_n$ be the tree whose nodes are the permutations of $S_n$ such that:
\begin{itemize}
\item $T_n$ has root $id_n$;
\item for every $\sigma, \tau\in S_n$, $\tau$ is a child of $\sigma$ if and only if $\bub(\tau)=\sigma$ and $\sigma\neq\tau$. \\

(Note that the situation $\bub(\tau)=\sigma$ and $\sigma = \tau$ occurs only when $\sigma = \tau = id_n$.)
\end{itemize}

Also, given a permutation $\pi\in S_n$, we define the \emph{tree of its preimages} $T(\pi)$ as the subtree of $T_n$ with root $\pi$.
\end{definition}

\begin{figure}[ht]
\centering
\begin{tikzpicture}
\tikzstyle{level 1}=[sibling distance=20mm]
\tikzstyle{level 2}=[sibling distance=10mm]
\tikzstyle{level 3}=[sibling distance=10mm]
\node {1234}
	child {node {1243}
	}
	child {node {1324}
		child {node {1342}
		}
		child {node {1432}
		}
		child {node {3142}
		}
		child {node {4132}
		}
	}
	child {node {1423}
	}
	child {node {2134}
		child {node {2314}
			child {node {4231}
			}
			child {node {2431}
			}
			child {node {2341}
			}
			child {node {3241}
			}
		}
		child {node {2413}
		}
		child {node {4213}
		}
		child {node {3214}
			child {node {4321}
			}
			child {node {3421}
			}
		}
	}
	child {node {2143}
	}
	child {node {3124}
		child {node {3412}
		}
		child {node {4312}
		}
	}
	child {node {4123}
	}
;
\end{tikzpicture}
\caption{The tree $T_4$.}
\label{T_4}
\end{figure}
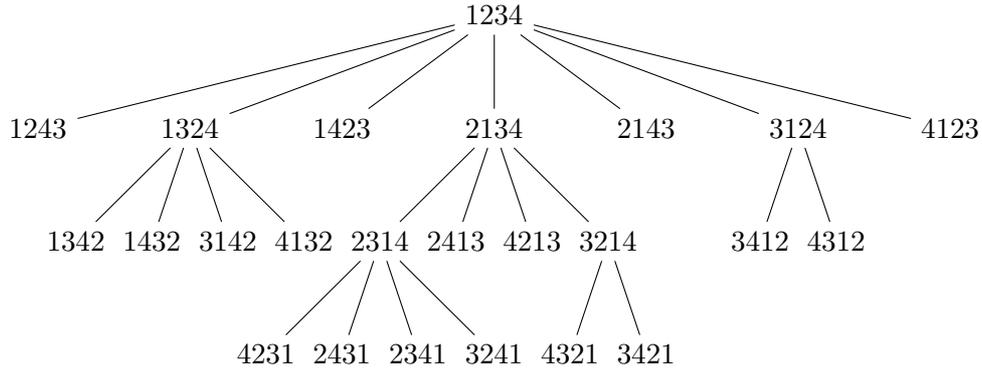

For example, \cref{T_4} shows the tree $T_4$.

\subsection{Isomorphisms between subtrees}

For any permutation $\pi$, $T(\pi)$ describes all possible preimages of $\pi$ under repeated applications of $\bub^{-1}$. 
From \cref{sec:CompPreimages}, we can see that the ``shape'' of this tree depends on $\pi$ only through the location of its left-to-right maxima. 
More precisely, the following lemma holds. 

\begin{lemma}\label{lem:same_LRmax_same_tree}
Let $\pi$ and $\tau$ be two permutations of the same size. 
If $\pi$ and $\tau$ have their left-to-right maxima in the same positions, then  $T(\pi)$ and $T(\tau)$ are isomorphic.
\end{lemma}
\begin{proof}
Let $\pi$, $\tau$ be permutations of the same size $n$ which have their left-to-right maxima in the same positions, and let $h$ and $k$ be the depth of $T(\pi)$ and $T(\tau)$, respectively (the \emph{depth} of a tree being defined as the maximum depth of its nodes). Without loss of generality, we can assume that $h\ge k$. The proof is by induction on $h$. If $h=0$, then $h=k=0$, and so $T(\pi)$ and $T(\tau)$ both consist of a single node, and our claim trivially holds.

Now suppose that $h>0$. 
Unless $\pi = id_n$, 
by definition, the nodes of $T(\pi)$ at depth $1$ are the preimages of $\pi$ under $\bub$. 
By \cref{cor:numberOfPreimages}, these preimages are in bijection with all possible subsets of left-to-right maxima of $\pi$ which do not contain $n$. 
Instead of identifying a subset of left-to-right maxima of $\pi$ by the \emph{values} of the left-to-right maxima it contains, 
we can identify it by the \emph{positions} of the left-to-right maxima it contains. 
Since $\pi$ and $\tau$ have their left-to-right maxima in the same positions, 
it follows from \cref{cor:numberOfPreimages} that there is a bijection between the preimages of $\pi$ under $\bub$ and  the preimages of $\tau$ under $\bub$. 
In addition, for every $\sigma \in \bub^{-1}(\pi)$, the corresponding $\rho \in \bub^{-1}(\tau)$ has its left-to-right maxima in the same positions as those of $\sigma$. 
Since $T(\sigma)$ and $T(\rho)$ have depth at most $h-1$ and $k-1$ respectively, 
we can apply the inductive hypothesis and obtain that $T(\sigma)$ and $T(\rho)$ are isomorphic. 
Summing up, we have that the children of $\pi$ and $\tau$ are in a bijective correspondence, 
and the trees rooted at two children paired together by this bijection are isomorphic. Therefore $T(\pi)$ and $T(\tau)$ are isomorphic.

Finally, if $\pi$ is the identity permutation of size $n$, 
and $\tau$ (of the same size) has its left-to-right maxima in the same positions as $\pi$, 
then necessarily $\tau = id_n$ as well, and in this case our claim trivially holds. 
\end{proof}

As the next proposition shows, all possible shapes of the trees $T(\pi)$ can be found starting at depth $1$ in $T_n$. 
To establish this proposition, we rely on the following decomposition of $\pi$, which will be also essential in the description of $T(\pi)$ in the next subsection. 

\begin{definition}
\label{dfn:decomposition}
Given a permutation $\pi$, we decompose it as $\pi=M_1 P_1 M_2 P_2 \cdots M_{\ell-1} P_{\ell-1} M_\ell$, 
where the $M_i$'s are all the maximal sequences of consecutive left-to-right maxima of $\pi$ (called \emph{blocks}), and the $P_i$'s collect all the remaining elements.
In particular, all the $P_i$'s are nonempty, and $M_i$ is nonempty for all $i$ except possibly for $i=\ell$.
Moreover, $m_i =|M_i |$ denotes the length of $M_i$, and analogously $p_i =|P_i |$ denotes the length of $P_i$, for all $i$.
\end{definition}

Notice that $m_1+\cdots +m_\ell=k$ is the total number of left-to-right maxima of $\pi$.

\begin{prop}\label{prop:isomorphicTrees}
For every permutation $\pi\in S_n$, $\pi\neq id_n$, there exists a child $\tau$ of $id_n$ in $T_n$ such that $T(\pi)$ and $T(\tau)$ are isomorphic. 
\end{prop}
\begin{figure}[ht]
\centering
\footnotesize
\begin{tikzpicture}[dot/.style={fill=black,circle}, scale=0.5]
\draw (0,0) rectangle (10,10);
\draw (0,0) -- (0.95,0.95);
\draw[loosely dotted] (1,-1) -- (1,3) node[shape=circle,fill=black,inner sep=1pt] {};
\draw[loosely dotted] (0,3) -- (3,3);
\draw (1.05,1) -- (3,2.95);
\draw[loosely dotted] (3,0) -- (3,3);
\draw (3.05,3.05) -- (5.95,5.95);
\draw[loosely dotted] (6,-1) -- (6,8) node[shape=circle,fill=black,inner sep=1pt] {};
\draw[loosely dotted] (0,8) -- (8,8);
\draw (6.05,6) -- (8,7.95);
\draw[loosely dotted] (8,0) -- (8,8);
\draw[loosely dotted] (8.05,8.05) -- (9.95,9.95);
\node[anchor=north] at (1,-1) {$m_1$};
\node[anchor=north] at (3,0) {$m_1+p_1$};
\node[anchor=north] at (6,-1) {$m_1+p_1+m_2$};
\node[anchor=north west] at (7.2,0) {$m_1+p_1+m_2+p_2$};
\node[anchor=east] at (0,3) {$m_1+p_1$};
\node[anchor=east] at (0,8) {$m_1+p_1+m_2+p_2$};
\end{tikzpicture}
\caption{The permutation $\tau$ described in the proof of \cref{prop:isomorphicTrees}.}
\label{tauExample}
\end{figure}
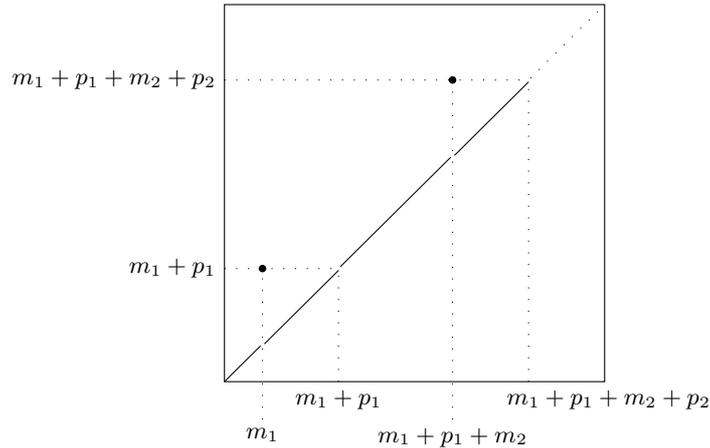
\begin{proof}
Let $\pi=M_1 P_1\cdots M_{\ell-1} P_{\ell-1} M_\ell \in S_n$. Define $\tau$ as the permutation in which the elements $m_1+p_1$, $m_1+p_1+m_2+p_2$, \dots , $m_1+p_1+\dots + m_{\ell-1} + p_{\ell-1}$, are in the positions $m_1$, $m_1+p_1+m_2$,\dots , $m_1+p_1+m_2+\dots + m_{\ell-2} + p_{\ell-2} + m_{\ell-1}$, respectively, while all the other elements are in increasing order. We can see an example of this construction in \cref{tauExample}. Therefore, $\tau$ and $\pi$ have their left-to-right maxima in the same positions, thus by \cref{lem:same_LRmax_same_tree} the trees $T(\pi)$ and $T(\tau)$ are isomorphic.

We are left with showing that $\tau$ is a child of $id_n$ in $T_n$. Since $p_1\neq 0$, then $\tau\neq id_n$, so we only need to check that $\bub(\tau) = id_n$. 
Observe that the elements $m_1+p_1+\dots +m_i+p_i$ are the last left-to-right maxima of their blocks in $\tau$, for every $i=1,\dots , \ell-1$, and all the elements before the positions $m_1+p_1+\dots +m_i+p_i$ are smaller than or equal to $m_1+p_1+\dots +m_i+p_i$. Therefore $\bub(\tau)=id_n$, because the $m_1+p_1+\dots +m_i+p_i$'s are the only elements moved by bubblesort, and they are moved to their correct position.
\end{proof}

\subsection{The skeleton of the tree of preimages}

Here we describe how the ``shape'' of any tree $T(\pi)$ is completely determined by a small piece of information about $\pi$, which we encapsulate in its \emph{label}. 

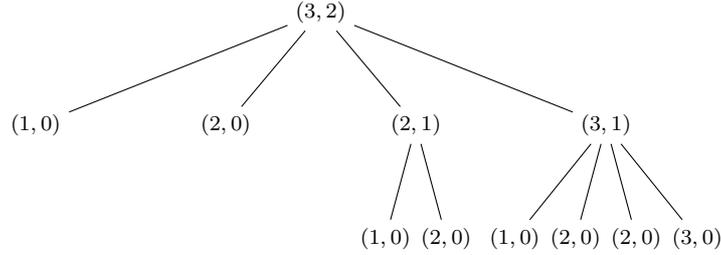
\begin{figure}[ht]
\centering
\scriptsize
\begin{tikzpicture}
\tikzstyle{level 1}=[sibling distance=25mm]
\tikzstyle{level 2}=[sibling distance=8mm]
\node {$(3,2)$}
	child {node {$(1,0)$}
	}
	child {node {$(2,0)$}
	}
	child {node {$(2,1)$}
		child {node {$(1,0)$}
		}
		child {node {$(2,0)$}
		}
	}
	child {node {$(3,1)$}
		child {node {$(1,0)$}
		}
		child {node {$(2,0)$}
		}
		child {node {$(2,0)$}
		}
		child {node {$(3,0)$}
		}
	}
;
\end{tikzpicture}
\caption{The skeleton of the tree $T(2134)$.}
\label{T(pi)}
\end{figure}

\begin{definition} 
The \emph{label} of a permutation $\sigma$ is the pair $(k, m_\ell)$, where $k$ and $m_\ell$ are defined as in \cref{dfn:decomposition}.\footnote{In particular, by definition of $k$ and $m_\ell$, the first component of a label is always at least as large as the second, with equality only in the case of the identity permutations.}
The \emph{skeleton} of a tree $T(\pi)$ is obtained from $T(\pi)$ by replacing each permutation at a node with its label. \cref{T(pi)} shows the skeleton of the tree $T(2134)$, and can be compared with the subtree $T(2134)$ of $T(1234)$ in \cref{T_4}.
\end{definition}

\begin{oss}
Since they have their left-to-right maxima in the same positions, the permutations $\pi$ and $\tau$ of \cref{prop:isomorphicTrees} have the same label. 
(It can also be observed that the trees $T(\pi)$ and $T(\tau)$ have the same skeleton. This follows by recursively applying \cref{cor:numberOfPreimages}, as in the proof of \cref{lem:same_LRmax_same_tree}.)
\end{oss}

Given a permutation $\pi$, we can determine the skeleton of $T(\pi)$ using only the pair $(k,m_\ell)$. 
Specifically, it is the tree with root labeled  by $(k,m_\ell)$, and whose children (and recursively, descendants) are obtained as described in the next proposition.

\begin{prop}
\label{prop:treegen}
Let $\pi\in S_n$ with label $(k,m_\ell)$. 
Let $T$ be the skeleton of $T(\pi)$.
Then the root of $T$ has label $(k,m_\ell)$ and its children have the following labels:
\begin{itemize}
\item for every $h=0,\ldots ,m_\ell -2$:
\begin{itemize}
   \item for every $i=1,\ldots ,k-1-h$, there are $\binom{k-2-h}{i-1}$ children with label $(k-i,h)$;
\end{itemize}
\item if $\pi\neq id_n$, we also have the case corresponding to $h=m_\ell -1$: 
\begin{itemize}
\item for every $i=0,\ldots ,k-m_\ell$, there are $\binom{k-m_\ell}{i}$ children with label $(k-i,m_\ell -1) = (k-i,h)$.
\end{itemize}
\end{itemize}
\end{prop}

\begin{proof}
We want to find the number of preimages of $\pi$ with any given label. 
If $m_\ell=0$, then $\pi$ does not end with its maximum, hence it has no preimage. Thus the root of $T$ has no children, and our claim vacuously holds.

Suppose that $m_\ell>0$. This means that $\pi=\pi_1\cdots \pi_{n - m_\ell} (n-m_\ell+1) \cdots n$. We can apply the procedure described in \cref{sec:CompPreimages} to find the preimages of $\pi$. 
From this procedure we can see that, if $\pi\neq id_n$, then its preimages can only have labels $(k',h)$ with $k' \leq k$ and $h <m_\ell$, which corresponds to the labels listed in the above statement.

If instead $\pi= id_n$, then it has label $(n,n)$ and its preimages can only have labels $(k',h)$, with $k' \leq n$ and $h \le n$, $h\neq n -1$. Indeed, we cannot obtain a preimage of $id_n$ with $h=n-1$, because that would mean that only the element 1 is not part of the last sequence of left-to-right maxima, which is impossible. Instead we can have $h=n$, but only by leaving $id_n$ unchanged.

To obtain a permutation with label $(k-i,h)$ with $h<m_\ell-1$, referring to the procedure of \cref{sec:CompPreimages}, we are forced to leave unchanged (i.e. not to swap) all the elements from $n$ down to $n-h+1$, then to swap $n-h$ with $n-h-1$. We are allowed to do so, because they are all left-to-right maxima, since $h<m_\ell-1$. 

After these steps, we obtain $\pi' (n-h-1) (n-h+1) (n-h+2) \cdots n$, with\footnote{Note that $\pi'$ is not a permutation, but just a sequence of distinct integers. However, as we have seen in \cref{lem:image}, it still makes sense to consider $\bub$ on such sequences.} $\pi'=\pi_1\cdots\pi_{n-h-2} (n-h)$.
Note that $\pi'$ has $k-h-1$ left-to-right maxima. 
Moreover, there is a bijection between the preimages of $\pi'$ 
and the preimages of $\pi$ ending with the suffix $(n-h-1) (n-h+1) (n-h+2) \cdots n$, which consists of just appending the suffix $(n-h-1) (n-h+1) (n-h+2)$. 
Under this bijection, if a preimage of $\pi'$ has $k-i-h$ left-to-right maxima (for some $i$ such that $1 \leq i \leq k-1-h$), 
then the corresponding preimage of $\pi$ has label $(k-i,h)$. 

From \cref{cor:PreimagesFixedNumberLTR}, the number of preimages of $\pi'$ with $k-i-h$ left-to-right maxima is $\binom{k-h-2}{k-i-h-1}=\binom{k-h-2}{i-1}$, for every $h=0,\ldots ,m_\ell -2$ and $i=1,\ldots ,k-1-h$. 
Exploiting the above bijection, this prove the first item of our proposition.

Consider now the case $h=m_\ell-1$, $\pi\neq id_n$. 
Then, applying the procedure of \cref{sec:CompPreimages}, 
we are forced to leave unchanged all the elements from $n$ down to $n-m_\ell +2$, then swap $n-m_\ell +1$ with $\pi_{n-m_\ell}$. Note that $\pi_{n-m_\ell}$ is an element of $\pi$ which is not a left-to-right maximum, so if we define $\pi'=\pi_1\cdots\pi_{n-m_\ell-1} (n-m_\ell+1)$, we have that $\pi'$ has $k-m_\ell+1 = k-h$ left-to-right maxima. An argument analogous to the one we have used for the case $h<m_\ell-1$ shows that the number of preimages of $\pi$ with label $(k-i, m_\ell -1)$ 
are $\binom{k-m_\ell}{k-i-m_\ell}=\binom{k-m_\ell}{i}$.

Finally, note that, if $\pi=id_n$, then there is an additional preimage, which is $id_n$, with label $(n,n)$. However it does not correspond to a child of $id_n$, because \cref{def:tree} prevents a permutation from being a child of itself.
\end{proof}

\begin{cor}
\label{cor:depth}
Let $\pi$ be a permutation of size $n$ with label $(k,m_\ell)$ such that $\pi\neq id_n$. Then $T(\pi)$ has depth $m_\ell$.
In addition, for every $n\ge1$, $T_n$ has depth $n -1$.
\end{cor}
\begin{proof}
We prove the statement by induction on $m_\ell$. 
If $m_\ell=0$, then $T(\pi)$ consists only of the root, and so has depth 0, as required.
If $m_\ell \ge 1$, then by \cref{prop:treegen} the root of $T(\pi)$ has children whose labels are of the form $(k',h)$, for every $0 \leq h\leq m_\ell -1$ and some $k'$. By induction hypothesis, the subtree rooted at each child with label $(k',h)$ has depth $h$. It follows that $T(\pi)$ has depth $1 + (m_\ell - 1) = m_\ell$, since the maximum value of $h$ is $m_\ell -1$.

We now consider $T_n$. If $n=1$ then $T_n$ consists of a single node and the statement is true. Otherwise, if $n>1$, then (again by \cref{prop:treegen}) the root of $T_n$ has children with labels $(k',h)$ for every $0 \leq h\le n-2$ and some $k'$. Since these children are not the identity permutation, we can apply the first part of this corollary to them. We obtain that each child with label $(h,k')$ is the root of a subtree of depth $h$. Since the maximum value of $h$ is $n-2$, $T_n$ has depth $1 + (n-2) = n-1$.
\end{proof}

\begin{cor}
For any given node $\pi\neq id_n$ in $T_n$, either half of its children are leaves or all of its children are leaves. 
\end{cor}
\begin{proof}
Let $(k,m_\ell)$ be the label of $\pi\neq id_n$. If $m_\ell=0$, the statement is vacuously true, because $\pi$ is a leaf.
Otherwise, if $m_\ell=1$, then by \cref{prop:treegen} we have that all of its children have label $(k',0)$ for some $k'$, and so they are all leaves.

Finally, if $m_\ell>1$, then (again by \cref{prop:treegen}) the number of children of $\pi$ which are leaves, that is with label $(k',0)$ for some $k'$, is
\[
\sum_{i=1}^{k-1} \binom{k-2}{i-1} = \sum_{j=0}^{k-2}\binom{k-2}{j}=2^{k-2}.
\]

By \cref{cor:PreimagesFixedNumberLTR}, $\pi$ has $2^{k-1}$ preimages, or equivalently it  has $2^{k-1}$ children (since $\pi\neq id_n$). This proves our statement.
\end{proof}

Notice that, for $\pi=id_n$ (with label $(n,n)$), it is still true that it has $2^{n-2}$ children which are leaves and $2^{n-1}$ preimages, but the total number of children is now $2^{n-1}-1$ ($id_n$ being a preimage of itself, but not one of its children).

\subsection{The inverse problem: deciding if a tree is isomorphic to $T(\pi)$ for some $\pi$}

We now consider the following problem: given a (rooted unlabeled) tree $T$, does there exist a permutation $\pi$ such that $T$ coincides with the (unlabeled) skeleton of $T(\pi)$? 
This problem can be easily solved thanks to our previous results on the labels of the nodes of $T(\pi)$. 
If $T$ consists of just a leaf, then of course $T$ is isomorphic to $T(\pi)$ for some $\pi$ (just take $\pi=1$ or any permutation of size at least $2$ not ending with its maximum). 
So, assume that $T$ has depth at least $1$. 

The first step is to determine the label of a candidate $\pi$. By \cref{cor:PreimagesFixedNumberLTR}, we can immediately say that, if the root of $T$ has neither $2^{k-1}$ nor $2^{k-1}-1$ children, for some $k>0$, then $T$ cannot be the (unlabeled) skeleton of any $T(\pi)$. 

If the root of $T$ has $2^{k-1}-1$ children, then the only candidate permutation is $\pi=id_k$. 
In particular, by \cref{cor:depth}, it is necessary that $T$ has depth $k-1$. 
We can then use \cref{prop:treegen} to check if the number of children of every node of $T$ matches with the numbers given in that proposition. 

Otherwise, suppose that the number of children of the root of $T$ is $2^{k-1}$ for some $k$. Let $m_\ell$ be the depth of $T$. 
By \cref{cor:PreimagesFixedNumberLTR} and \cref{cor:depth}, we know that a permutation needs to have label $(k,m_\ell)$ 
for $T(\pi)$ to have $2^{k-1}$ children of its root and depth $m_\ell$.
Therefore we can use \cref{prop:treegen} to check if the number of children of every node of $T$ matches with the numbers given in \cref{prop:treegen} for a permutation $\pi\neq id_n$ with label $(k,m_\ell)$.

The next proposition summarizes the above discussion.

\begin{prop}
Let $T$ be a (rooted unlabeled) tree, let $i$ be the number of children of the root of $T$, and $m_\ell$ be the depth of $T$. Then
\begin{itemize}
\item if $i=2^{m_\ell}-1$, then $T$ may only coincide with the (unlabeled) skeleton of $T_{m_\ell +1}$;
\item if there exists a positive integer $k$ such that $i=2^{k-1}$, then $T$ may only coincide with the (unlabeled) skeleton of a permutation $\pi$ with label $(k,m_\ell)$;
\item in all the other cases, $T$ does not coincide with the unlabeled skeleton of any permutation.
\end{itemize}
\end{prop}

\section{Heights of nodes and leaves in $T_n$} 
\label{sec:height}

\subsection{Nodes}

Recall that the height of a node in a rooted tree is the number of edges on the path connecting that node to the root. 
The height of a node of the tree $T_n$ corresponds to the number of passes of Bubblesort needed to sort the permutation at this node.
Therefore, we can refer to \cite[Prop. 17]{bubbleM} to find information on the number of nodes of $T_n$.

\begin{prop}[\cite{bubbleM}]
\label{prop:AABCD}
The set of permutations of size $n$ sorted by at most $k$ passes of Bubblesort is the set $\textnormal{Av}_n(\Gamma_{k+2})$, where $\Gamma_k$ is the set of all permutations of size $k$ whose final element is 1\footnote{We warn the reader that we have made a slight change of notation with respect to \cite{bubbleM} here; more specifically, our set $\Gamma_k$ is $\Gamma_{k-2}$ in \cite{bubbleM}.}. As a consequence, setting $\varphi_n^{(k)}=|\textnormal{Av}_n(\Gamma_{k})|$, the number of nodes at height at most $k$ in $T_n$ is given by $\varphi_{n}^{(k+2)}=(k+1)^{n-k-1}(k+1)!$
\end{prop}

We can thus immediately deduce the number of nodes at a given height in $T_n$.

\begin{cor}
\label{cor:AABCD}
The number $f_n^{(k)}$ of nodes at height $k$ in $T_n$ is given by	
\[
f_n^{(k)}=\varphi_{n}^{(k+2)}-\varphi_{n}^{(k+1)}=(k+1)^{n-k-1}(k+1)!-k^{n-k}k!= k!\cdot ((k+1)^{n-k}-k^{n-k}).
\]
\end{cor}

The first lines of the infinite triangular matrix of the coefficients $f_n^{(k)}$ are given in Table \ref{nodes_height}. This is sequence A056151 in \cite{Sl}.

\begin{table}
\begin{center}
	\begin{tabular}{c|cccccc}
		\backslashbox{$n$}{$k$} & 0 & 1 &2&3&4&5\\
		\hline
		1 & 1 &   &   &   &    &          \\
		2 & 1 & 1 &   &   &    &          \\
		3 & 1 & 3 & 2 &   &    &          \\
		4 & 1 & 7 & 10 & 6 &    &         \\
		5 & 1 & 15 & 38 & 42 & 24    &    \\
		6 & 1 & 31 & 130 & 222 & 216 & 120    \\
	\end{tabular}
	\end{center}
	\caption{Number of nodes in $T_n$ having height $k$.} \label{nodes_height}
\end{table}

We notice that the elements on the diagonal of Table \ref{nodes_height} are the factorial numbers, more specifically $f_{n}^{(n-1)}=(n-1)!$. Indeed, the set of permutations of size $n$ needing the maximum number of passes of Bubblesort to be sorted (that is, $n-1$ passes) is the set of permutations of size $n$ ending with 1, whose cardinality is clearly $(n-1)!$.

\medskip

From the expression of $\varphi_n^{(k)}$ in \cref{prop:AABCD}, 
we can derive the asymptotic behavior of the average height of a node in $T_n$. 
This analysis is described in~\cite[Theorem 7.14]{FS} 
and follows easily from the asymptotic behavior of the Ramanujan P-function (see~\cite[Table 4.11]{FS} or~\cite[p. 119-120]{knuth1}), 
which we state in \cref{lem:Pfunction} below. 
We then reproduce the analysis of~\cite[Theorem 7.14]{FS}, as a preparation for \cref{prop:average_height_leaves} below. 

\begin{lemma}[\cite{FS}]
\label{lem:Pfunction}
The Ramanujan P-function, defined by $P(n) = \sum_{k=0}^{n-1} \frac{k! k ^{n-k}}{n!}$, behaves asymptotically as $\sqrt{\tfrac{\pi n}{2}} +O(1)$.
\end{lemma}

\begin{prop}[\cite{FS}]
\label{prop:average_height_nodes}
The average height of a node in $T_n$ is asymptotically equal to $n-\sqrt{\frac{\pi n}{2}} + O(1)$.
\end{prop}

\begin{proof}
The average height of a node in $T_n$ is given by 
\[
H_n := \frac{1}{n!}\sum_{k=1}^{n-1} \text{number of nodes of height at least } k \text{ in } T_n,
\]
each node at height $k$ contributing indeed exactly $k$ times to this sum. 
Writing the number of nodes of height at least $k$ in $T_n$ as the difference of $n!$ (the total number of nodes) and the number of nodes of height at most $k-1$ in $T_n$, we then compute 
\[
 H_n = \frac{1}{n!}\sum_{k=1}^{n-1} (n! - \varphi_n^{(k+1)})
 = \frac{1}{n!}\sum_{k=1}^{n-1} (n! - k^{n-k}k!) 
 = (n-1) - \sum_{k=0}^{n-1} \frac{k^{n-k}k!}{n!}
 = (n-1) - P(n),
\]
proving our claim. 
\end{proof}

\bigskip

Recall the (obvious) fact that $T_n$ contains $n!$ nodes. 
With \cref{prop:AABCD,cor:AABCD}, we have refined this counting according to the height of the nodes in $T_n$. 
We now address the analogous problems in $T(\pi)$ for $\pi \neq id_n$. 
More precisely, given a permutation $\pi$ (of size $n$) having label $(k,m_\ell)$, we determine an expression for the number of nodes of $T(\pi)$ (which does not depend on $n$). 
This expression is a summation formula in which each summand counts nodes in $T(\pi)$ of a prescribed height.

\begin{lemma}\label{lem:nodes_level}
	Let $\pi$ and $\tau$ be two permutations having labels $(k,m_\ell)$ and $(k,m_\ell -1)$, respectively, with $1\leq m_\ell \leq k-1$. 
	Then the tree obtained by removing the leaves at height $m_\ell$ in $T(\pi)$ is isomorphic to $T(\tau)$.
\end{lemma}

\begin{proof}
	Remember that, by \cref{cor:depth}, $T(\pi)$ has height $m_\ell$ and $T(\tau)$ has height $m_\ell-1$. The proof is by induction on $m_\ell$.
	
	If $m_\ell=1$, then $T(\tau)$ consists of the single node $\tau$, while $T(\pi)$ has height $1$, therefore the statement is true.
	
	Now let $m_\ell\ge 2$, and suppose that the statement is true for $m_\ell-1$. We will show that there is a bijective correspondence between the children of $\tau$ and the children of $\pi$ such that the subtree rooted at a child of $\tau$ is isomorphic to the subtree rooted at the corresponding child of $\pi$, after removing the leaves at height $m_\ell$ (if any). 
	
	\cref{prop:treegen} allows us to determine the labels of the children of $\tau$ and $\pi$ in $T(\tau)$ and $T(\pi)$, respectively. Specifically, $\tau$ and $\pi$ have the same number of children with labels $(k-i,h)$, for every $h=0,\dots,m_\ell-3$ and every $i=0,\dots,k-1-h$. Regarding the remaining children, we have that the number of children of $\tau$ labeled $(k-i,m_\ell-2)$ is equal to the sum of the number of children of $\pi$ labeled $(k-i,m_\ell-2)$ and $(k-i,m_\ell-1)$, for every $i=0,\ldots ,k-m_\ell+1$.
This induces the announced bijective correspondence between the children of $\tau$ in $T(\tau)$ and those of $\pi$ in $T(\pi)$. 
	
	The children of $\tau$ and $\pi$ with the same labels give isomorphic subtrees by \cref{prop:treegen}. 
	In addition, if this label is $(k-i,h)$ for some $h$ $\leq$ $m_\ell -2$ (and some suitable $i$), then the subtrees contain no leaf at height $m_\ell$ in $T(\tau)$ or $T(\pi)$ (again by \cref{cor:depth}), ensuring our claim restricted to such children of $\pi$ and $\tau$. 

	Therefore, we are left with considering a child of $\pi$ in $T(\pi)$ with label $(k-i,m_\ell-1)$, to which corresponds a child of $\tau$ in $T(\tau)$ of label  $(k-i,m_\ell-2)$. 
	We can apply the inductive hypothesis to such children of $\pi$ and $\tau$, thus obtaining that each subtree of $T(\tau)$ rooted at a child of $\tau$ with label $(k-i,m_\ell-2)$ is isomorphic to a subtree of $T(\pi)$ rooted at a child of $\pi$ with label $(k-i,m_\ell-1)$ after removing the leaves at height $m_\ell-1$ (in the subtree, \emph{i.e.} at height $m_\ell$ in $T(\pi)$).
This concludes the proof.
\end{proof}

\begin{prop}
\label{nodes_T(pi)}
For a permutation $\pi$ having label $(k,m_\ell)$, different from an identity permutation, 
the number of nodes of the tree $T(\pi)$ of its preimages under $\bub$ is
\begin{equation}\label{nodes}
N(k,m_\ell )=\sum_{j=0}^{m_\ell}j!(j+1)^{k-j}.
\end{equation}

Moreover, each summand in \cref{nodes} records the contribution of each level of $T(\pi )$. In other words, denoting with $N_j (k,m_\ell )$ the number of nodes at height $j$ in $T(\pi)$, we have that $N_j (k,m_\ell )=j!(j+1)^{k-j}$.
\end{prop}

\begin{proof}
In order to prove \cref{nodes} we proceed by induction on $m_\ell$. If $m_\ell =0$, then $\pi$ has no children, hence $N(k,0)=1$, which is consistent with \cref{nodes}.

Now suppose that \cref{nodes} holds when the cardinality of the longest suffix of left-to-right maxima of $\pi$ is strictly smaller than $m_\ell$. Recalling \cref{prop:treegen}, we have the following recursive expression for the number of nodes of $T(\pi)$:

{\small
\begin{align*}
N(k,m_\ell )&=1+\sum_{h=0}^{m_\ell -2}\sum_{i=1}^{k-1-h}\binom{k-2-h}{i-1}N(k-i,h)+\sum_{i=0}^{k-m_\ell}\binom{k-m_\ell}{i}N(k-i,m_\ell -1) \\
&=1+\sum_{h=0}^{m_\ell -2}\sum_{i=1}^{k-1-h}\binom{k-2-h}{i-1}\sum_{j=0}^{h}j!(j+1)^{k-i-j}+\sum_{i=0}^{k-m_\ell}\binom{k-m_\ell}{i}\sum_{j=0}^{m_\ell -1}j!(j+1)^{k-i-j} \\
&=1+\sum_{h=0}^{m_\ell -2}\sum_{j=0}^{h}j!(j+1)^{k-j-1}\sum_{i=0}^{k-2-h}\binom{k-2-h}{i}(j+1)^{-i}+\sum_{j=0}^{m_\ell -1}j!(j+1)^{k-j}\sum_{i=0}^{k-m_\ell}\binom{k-m_\ell}{i}(j+1)^{-i} \\
&=1+\sum_{h=0}^{m_\ell -2}\sum_{j=0}^{h}j!(j+1)^{k-j-1}\left( 1+\frac{1}{j+1}\right)^{k-2-h}+\sum_{j=0}^{m_\ell -1}j!(j+1)^{k-j}\left( 1+\frac{1}{j+1}\right) ^{k-m_\ell} \\
&=1+\sum_{h=0}^{m_\ell -2}\sum_{j=0}^{h}j!(j+1)^{h-j+1}(j+2)^{k-2-h}+\sum_{j=0}^{m_\ell -1}j!(j+1)^{m_\ell -j}(j+2)^{k-m_\ell}.
\end{align*}
}

We then exchange the order of the two sums in the middle term of the last expression, use the geometric sum formula and we get:

\begin{align*}
N(k,m_\ell )&=1+\sum_{j=0}^{m_\ell -2}j!(j+1)^{1-j}(j+2)^{k-2}\sum_{h=j}^{m_\ell -2}(j+1)^{h}(j+2)^{-h}+\sum_{j=0}^{m_\ell -1}j!(j+1)^{m_\ell -j}(j+2)^{k-m_\ell} \\
&=1+\sum_{j=0}^{m_\ell -2}j!(j+1)(j+2)^{k-1-j}-\sum_{j=0}^{m_\ell -2}j!(j+1)^{m_\ell -j}(j+2)^{k-m_\ell}+\sum_{j=0}^{m_\ell -1}j!(j+1)^{m_\ell -j}(j+2)^{k-m_\ell} \\
&=1+\sum_{j=1}^{m_\ell -1}j!(j+1)^{k-j}+m_\ell !(m_\ell +1)^{k-m_\ell}=\sum_{j=0}^{m_\ell}j!(j+1)^{k-j}, 
\end{align*} 
which gives \cref{nodes}.

Concerning the evaluation of $N_j (k,m_\ell )$, \cref{lem:nodes_level} implies that $N_j (k,m_\ell)=N_j (k,m_\ell -1)$, for all $j\leq m_\ell -1$. By a repeated application of the lemma, we get that $N_j (k,m_\ell)=N_j (k,j)=N(k,j)-N(k,j-1)=j!(j+1)^{k-j}$, as desired.  
\end{proof}

\subsection{Leaves}

In the tree $T_n$ the leaves represent permutations that cannot be obtained as output of Bubblesort, \emph{i.e.}, which do not belong to the image of $\bub$. 
We saw just after \cref{lem:image} that these permutations are those not ending with their maximum, so that 
the total number of leaves in $T_n$ is given by $(n-1)\cdot (n-1)!$.  

Our next result is a closed formula for the number of leaves at height $k$ in $T_n$, for any $k\leq n-1$. To this aim, we make use of the so-called \emph{ECO method}, illustrated in \cite{BDLPP} and further developed and employed by many authors (see for instance \cite{FPPR}). We will not give a detailed description of this method here, since our application is simple enough to be outlined directly.

\bigskip

Recall that leaves in $T_n$ correspond to permutations whose last element is not the maximum. Thus, denoting with $\textnormal{Av}^*_n(\Gamma_{k})$ the set of permutations of size $n$ avoiding $\Gamma_k$ and such that their last element is different from $n$, we are interested in the coefficients $\gamma_n^{(k)}=|\textnormal{Av}^*_n(\Gamma_{k})|$, since $\gamma_n^{(k+2)}$ gives the number of leaves at height at most $k$ in $T_n$.

\begin{prop}
\label{prop:height_leaves}
	For all $n,k$, we have
	\[
	\gamma_n^{(k)}=
	\begin{cases}
	(n-1)(n-1)! & n<k, \\
	(k-2)(k-1)^{n-k}(k-1)! & n\geq k.
	\end{cases}
	\]
\end{prop}

\begin{proof}
	We consider the following general procedure to generate all permutations of size $n$. Given any permutation of size $n-1$, construct $n$ different permutations of size $n$ by adding a new rightmost element $k$, for any choice of $k$ between 1 and $n$, and suitably rescaling the other elements (namely, all elements of the starting permutation which are greater than or equal to $k$ are increased by 1, whereas all the remaining elements are left untouched). It is immediate to realize that, starting from the set of all permutations of size $n-1$, the above procedure generates exactly once every permutation of size $n$.
	
	We now adapt the above construction to our setting. Every permutation of $\textnormal{Av}^*_n(\Gamma_{k})$ can be obtained from a permutation of $\textnormal{Av}_{n-1}(\Gamma_{k})$ by adding a suitable rightmost element. More specifically, we cannot add $n$ (because we require that our permutation does not end with its maximum); moreover, if $n\geq k$, we cannot add any element between 1 and $n-k+1$ as well (otherwise we would create one of the forbidden patterns belonging to $\Gamma_k$). On the other hand, any of the remaining elements is allowed and generates a valid permutation. This means that every permutation in $\textnormal{Av}_{n-1}(\Gamma_{k})$ generates $k-2$ distinct permutations of $\textnormal{Av}^*_n(\Gamma_{k})$ and every permutation in $\textnormal{Av}^*_n(\Gamma_{k})$ is obtained in this way exactly once. We thus deduce that, when $n\geq k$,
	\[
	\gamma_n^{(k)}=(k-2)\varphi_{n-1}^{(k)}=(k-2)(k-1)^{n-k}(k-1)!,
	\]
	whereas for $n<k$ we have that $\gamma_n^{(k)}=(n-1)(n-1)!$, which concludes the proof.
\end{proof}

\begin{cor}
	The number $g_n^{(k)}$ of leaves of $T_n$ at height $k$ is given by
	\[
	g_n^{(k)}=k!(k(k+1)^{n-k-1}-(k-1)k^{n-k-1}).
	\]
\end{cor}

\begin{proof}
	Just observe that $g_n^{(k)}=\gamma_n^{(k+2)}-\gamma_n^{(k+1)}$ and that the maximum height of a node of $T_n$ is $n-1$, so we are only interested in the case $n\geq k+1$ of the previous proposition.
\end{proof}

As in the case of nodes, \cref{prop:height_leaves} allows us to derive the asymptotic behavior of the average height of a leaf in $T_n$. 

\begin{prop}
\label{prop:average_height_leaves}
The average height of a leaf in $T_n$ is asymptotically equal to $n-\sqrt{\frac{\pi n}{2}} + O(1)$.
\end{prop}

\begin{proof}
As in the proof of \cref{prop:height_leaves}, we have that 
the average height of a leaf in $T_n$ is  
\begin{align*}
G_n & = \frac{1}{(n-1)(n-1)!}\sum_{k=1}^{n-1} \text{number of leaves of height at least } k \text{ in } T_n \\
& = \frac{1}{(n-1) (n-1)!}\sum_{k=1}^{n-1} \big((n-1) (n-1)! - \gamma_n^{(k+1)}\big)
= (n-1) - \sum_{k=1}^{n-1} \frac{(k-1)k^{n-k-1}k!}{(n-1) (n-1)!} \\
& = (n-1) - \frac{n}{n-1} \sum_{k=1}^{n-1} \frac{k^{n-k}k!}{n!} 
+ \frac{1}{n-1} \sum_{k=1}^{n-1} \frac{k^{n-1-k}k!}{(n-1)!} 
= (n-1) - \frac{n P(n)}{n-1}  + \frac{P(n-1) +1}{n-1}, 
\end{align*}
and the asymptotic behavior of the Ramanujan P-function yields the announced result. 
\end{proof}

In the same manner as we have done for the nodes, we now address the analogous problem of counting the leaves in $T(\pi)$, for $\pi \neq id_n$.
More precisely, given a permutation $\pi$ (of size $n$) having label $(k,m_\ell)$, we determine an expression for the number of leaves of $T(\pi)$ (which does not depend on $n$ but only on the label $(k,m_\ell)$). 
This expression is a summation formula in which each summand counts the leaves of a prescribed height in $T(\pi)$.

\begin{prop}
\label{prop:leavesTpi}
For a permutation $\pi$ having label $(k,m_\ell)$, different from an identity permutation, the number of leaves of the tree $T(\pi)$ of its preimages under $\bub$ is
\begin{equation}\label{leaves}
L(k,m_\ell )=\sum_{j=1}^{m_\ell-1}j! j (j+1)^{k-j-1} + m_\ell !(m_\ell+1)^{k-m_\ell}.
\end{equation}

Moreover, each summand in \cref{leaves} records the contribution of each level of $T(\pi )$. In other words, denoting with $L_j (k,m_\ell )$ the number of leaves at height $j$ in $T(\pi)$, we have that $L_j (k,m_\ell )=j! j (j+1)^{k-j-1}$ for $j< m_\ell$, and $L_{m_\ell} (k,m_\ell ) = m_\ell !(m_\ell+1)^{k-m_\ell}$.
\end{prop}

\begin{proof}
The proof of \cref{leaves} is by induction, following the exact same steps as the proof of \cref{nodes_T(pi)}.
The recursive equation for the number of leaves in $T(\pi)$, which is needed in the inductive step of the proof, is again obtained from \cref{prop:treegen}. It actually differs from the one for nodes in the proof of \cref{nodes_T(pi)} only by the initial term $1$ (accounting for the root node); namely for $m_\ell \geq 1$, we have 
\[
L(k,m_\ell )=\sum_{h=0}^{m_\ell -2}\sum_{i=1}^{k-1-h}\binom{k-2-h}{i-1}L(k-i,h)+\sum_{i=0}^{k-m_\ell}\binom{k-m_\ell}{i}L(k-i,m_\ell -1),
\]
and for $m_\ell = 0$ it holds that $L(k,0)=1$.
From there, the same steps of computations as in the proof of \cref{nodes_T(pi)} (followed by additional elementary simplifications) yield, for $m_\ell \geq 1$:
\[
L(k,m_\ell )= \sum_{j=1}^{m_\ell -1}j! j (j+1)^{k-j-1}+m_\ell !(m_\ell +1)^{k-m_\ell},  
\]
as claimed. 

We now move to the claimed expression for $L_j (k,m_\ell )$. 
We shall first establish it for $j=m_\ell$, then for $j=m_\ell-1$, and then for smaller $j$ iterating the argument. 

We first note that all the nodes of $T(\pi)$ at height $m_\ell$ are leaves (since $m_\ell$ is the height of this tree). Using \cref{nodes_T(pi)}, we therefore have $L_{m_\ell} (k,m_\ell )= N_{m_\ell} (k,m_\ell ) = m_\ell !(m_\ell+1)^{k-m_\ell}$. As a consequence, the total number of leaves having height at most $m_\ell-1$ in $T(\pi)$ is $\sum_{j=1}^{m_\ell-1}j! j (j+1)^{k-j-1}$.

Next, we claim that the number of leaves having height at most $m_\ell-2$ in $T(\pi)$ is $\sum_{j=1}^{m_\ell-2}j! j (j+1)^{k-j-1}$. From this claim, the announced formula $L_{m_\ell-1}(k,m_\ell ) = (m_\ell-1)! (m_\ell-1) m_\ell^{k-m_\ell}$ immediately follows by taking the difference. 

To prove our claim, we use \cref{lem:nodes_level}. This lemma indeed implies that $L_j (k,m_\ell)=L_j (k,m_\ell-1)$, for all $j\leq m_\ell -2$. This shows that the number of leaves having height at most $m_\ell-2$ in $T(\pi)$ is the same as the number of leaves having height at most $m_\ell-2$ in $T(\sigma)$ for $\sigma$ a permutation with label $(k,m_\ell-1)$. The latter is equal to $L(k,m_\ell-1) - L_{m_\ell -1}(k,m_\ell-1)$, hence equal to $\sum_{j=1}^{m_\ell-2}j! j (j+1)^{k-j-1}$ as established earlier, thus proving our claim. 

We are now left with showing that $L_h (k,m_\ell )=h! h (h+1)^{k-h-1}$ for $h\leq m_\ell-2$. We proceed iteratively, for decreasing values of $h$. 
At each step, the reasoning is similar to the above case for $h=m_\ell -1$. 
We first use \cref{lem:nodes_level} (several times, as in the proof of \cref{nodes_T(pi)}) to argue that the number of leaves having height at most $h-1$ in $T(\pi)$ is the same as the number of leaves having height at most $h-1$ in $T(\sigma)$ for $\sigma$ a permutation with label $(k,h)$. This number is $\sum_{j=1}^{h-1}j! j (j+1)^{k-j-1}$.
Then, $L_h (k,m_\ell )$ is the difference between $L(k,m_\ell) - \sum_{h+1\leq j\leq m_\ell} L_j(k,m_\ell)$ and the above quantity. 
The result follows from the formulas previously established for $L_j(k,m_\ell)$ for $j \geq h+1$.
\end{proof}

\begin{remark}
Combining Propositions~\ref{nodes_T(pi)} and~\ref{prop:leavesTpi} tells us that, for $\pi$ a permutation of label $(k,m_\ell)$, at height $j<m_\ell$ in $T(\pi)$, the ratio between the number of leaves and the number of nodes is $\frac{j}{j+1}$ (equivalently, the ratio between the number of internal nodes and the number of nodes is $\frac{1}{j+1}$).
\end{remark}

\bibliographystyle{plain}
\bibliography{bubble}

\begin{thebibliography}{10}

\bibitem{bubbleM}
Michael~H. Albert, M.~D. Atkinson, Mathilde Bouvel, Anders Claesson, and Mark
  Dukes.
\newblock On the inverse image of pattern classes under bubble sort.
\newblock {\em Journal of Combinatorics}, 2(2):231--243, 2011.

\bibitem{BDLPP}
E.~Barcucci, A.~Del~Lungo, E.~Pergola, and R.~Pinzani.
\newblock {ECO}: a methodology for the enumeration of combinatorial objects.
\newblock {\em Journal of Difference Equations and Applications}, 5:435--490,
  1999.

\bibitem{BonaSurvey}
Mikl\'os B\'ona.
\newblock A survey of stack-sorting disciplines.
\newblock {\em The Electronic Journal of Combinatorics}, 9(2), 2003.

\bibitem{Mireille}
Mireille Bousquet-Mélou.
\newblock Sorted and/or sortable permutations.
\newblock {\em Discrete Mathematics}, 225(1--3):25--50, 2000.

\bibitem{LucaLapo2}
Lapo Cioni and Luca Ferrari.
\newblock Characterization and enumeration of preimages under the queuesort
  algorithm.
\newblock In J.~Nešetřil, G.~Perarnau, J.~Rué, and O.~Serra, editors, {\em
  Extended Abstracts EuroComb 2021. Trends in Mathematics}, volume~14 of {\em
  Birkhäuser, Cham}, 2021.

\bibitem{LucaLapo1}
Lapo Cioni and Luca Ferrari.
\newblock Preimages under the queuesort algorithm.
\newblock {\em Discrete Mathematics}, 344, 2021.

\bibitem{ClaessonUlfarsson}
Anders Claesson and Henning Ulfarsson.
\newblock Sorting and preimages of pattern classes.
\newblock {\em DMTCS Proceedings}, AR:595--606, 2012.

\bibitem{defant4}
Colin Defant.
\newblock Preimages under the stack-sorting algorithm.
\newblock {\em Graphs and Combinatorics}, 33:103--122, 2017.

\bibitem{defant3}
Colin Defant.
\newblock Fertility numbers.
\newblock {\em Journal of Combinatorics}, 11:527--548, 2020.

\bibitem{defant7}
Colin Defant.
\newblock Polyurethane toggles.
\newblock {\em Electronic Journal of Combinatorics}, 27(2):Article P2.46, 2020.

\bibitem{defant6}
Colin Defant.
\newblock Stack-sorting preimages of permutation classes.
\newblock {\em S\'eminaire Lotharingien de Combinatoire}, 82:Article B82b,
  2020.

\bibitem{defant5}
Colin Defant.
\newblock Enumeration of stack-sorting preimages via a decomposition lemma.
\newblock {\em Theoretical Computer Science}, 22(3), 2021.

\bibitem{defant1}
Colin Defant.
\newblock Fertility monotonicity and average complexity of the stack-sorting
  map.
\newblock {\em European Journal of Combinatorics}, 93, 2021.

\bibitem{defant2}
Colin Defant, Michael Engen, and Jordan~A. Miller.
\newblock Stack-sorting, set partitions, and {L}assalle's sequence.
\newblock {\em Journal of Combinatorial Theory, Series A}, 175, 2020.

\bibitem{FPPR}
L.~Ferrari, E.~Pergola, R.~Pinzani, and S.~Rinaldi.
\newblock Some applications arising from the interactions between the theory of
  {C}atalan-like numbers and the {ECO} method.
\newblock {\em Ars Combinatoria}, 99:109--128, 2011.

\bibitem{LucaTalk}
Luca Ferrari.
\newblock Sorting with stacks and queues: some recent developments.
\newblock Keynote address at the on-line conference \emph{Permutation Patterns
  2021}.
\newblock Available at \url{https://www.youtube.com/watch?v=cTT9t5gddmE}.

\bibitem{FS}
P.~Flajolet and R.~Sedgewick.
\newblock {\em An {I}ntroduction to the {A}nalysis of {A}lgorithms, second
  edition}.
\newblock Addison-Wesley, 2013.

\bibitem{knuth1}
Donald~E. Knuth.
\newblock {\em The Art of Computer Programming, Vol. 1: Fundamental
  Algorithms}.
\newblock Addison-Wesley, third edition, 1997.

\bibitem{Magnusson}
Hjalti Magnusson.
\newblock Sorting operators and their preimages.
\newblock Master's thesis, Reykjavik University, 2013.

\bibitem{Sl}
N.~J.~A. Sloane.
\newblock {T}he {O}nline {E}ncyclopedia of {I}nteger {S}equences.
\newblock \url{oeis.org}.

\bibitem{Tao}
Jiang Tao, Li~Ming, and Paul~M.B. Vitanyi.
\newblock Average-case analysis of algorithms using {K}olmogorov complexity.
\newblock {\em Journal of Computer Science and Technology}, 15:402--408, 2000.

\bibitem{VatterSurvey}
Vincent Vatter.
\newblock {\em Permutation Classes}, chapter 12 of {\em The Handbook of
  Enumerative Combinatorics}, pages 753--834.
\newblock Chapman and Hall/CRC Press, 2015.

\end{thebibliography}

\end{document}